\documentclass[10pt]{article}

\usepackage{amsmath}
\usepackage{amsfonts}
\usepackage{amssymb}

\newtheorem{proof}{Proof}

\newtheorem{theorem}{Theorem}[section]
\newtheorem{lemma}[theorem]{Lemma}
\newtheorem{corollary}[theorem]{Corollary}
\newtheorem{remark}[theorem]{Remark}

\newtheorem{definition}{Definition}

\def\phi{{\varphi}}

\DeclareSymbolFont{AMSb}{U}{msb}{m}{n}
\DeclareMathSymbol{\N}{\mathbin}{AMSb}{"4E}
\DeclareMathSymbol{\Z}{\mathbin}{AMSb}{"5A}
\DeclareMathSymbol{\R}{\mathbin}{AMSb}{"52}
\DeclareMathSymbol{\Q}{\mathbin}{AMSb}{"51}
\DeclareMathSymbol{\I}{\mathbin}{AMSb}{"49}
\DeclareMathSymbol{\C}{\mathbin}{AMSb}{"43}

\begin{document}
\title{On the critical dimension of a fourth order elliptic problem with negative exponent}
\author{  Amir Moradifam \thanks{This work is supported by a Killam Predoctoral
Fellowship, and is part of the author's PhD dissertation in
preparation under the supervision of N. Ghoussoub.}\\
\small Department of Mathematics,
\small University of British Columbia, \\
\small Vancouver BC Canada V6T 1Z2 \\
\small {\tt a.moradi@math.ubc.ca}
\\
} \maketitle

\begin{abstract}
We study the regularity of the extremal solution of the semilinear
biharmonic equation $\beta \Delta^2 u-\tau \Delta
u=\frac{\lambda}{(1-u)^2}$ on a ball $B \subset \R^N$, under Navier
boundary conditions $u=\Delta u=0$ on $\partial B$, where $\lambda
>0$ is a parameter, while $\tau>0$, $\beta>0$ are fixed constants.
It is known that there exists a $\lambda^{*}$ such that for
$\lambda>\lambda^{*}$ there is no solution while for
$\lambda<\lambda^{*}$ there is a branch of minimal solutions. Our
main result asserts that the extremal solution $u^{*}$ is regular
($\sup_{B}u^{*}<1$) for $N\leq 8$ and $\beta, \tau>0$ and it is
singular ($\sup_{B}u^{*}=1$) for $N\geq 9$, $\beta>0$, and $\tau>0$
with $\frac{\tau}{\beta}$ small. Our proof for the singularity of
extremal solutions in dimensions $N\geq 9$ is based on certain
improved Hardy-Rellich inequalities.
\end{abstract}

\section{Introduction}

Consider the fourth order elliptic problem
$$ \left \{ \begin{array}{cl}
\beta \Delta^2 u- \tau \Delta u=\frac{\lambda}{(1-u)^2} \;\;
&\mbox{in $\Omega$},\\ 0<u \leq 1\;\; &\mbox{in $\Omega$},\\
u=\Delta u=0 \;\; &\mbox{on $\partial \Omega$,}
\end{array} \right. \eqno(G_\lambda)$$ where $\lambda>0$ is a parameter,
$\tau>0$, $\beta>0$ are fixed constants, and $\Omega \subset \R^N\;
(N \geq 2)$ is a bounded  smooth domain. This problem with $\beta=0$
models a simple electrostatic Micro-Electromechanical Systems (MEMS)
device which has been recently studied by many authors.  For
instance, see \cite{CG}, \cite{CR}, \cite{EGG}, \cite{GaG},
\cite{GG1}, \cite{GG2}, \cite{G}, \cite{GW1}, \cite{GW2},
\cite{GW3}, and the references cited therein.

Recently, Lin and Yang \cite{LY} derived the equation
$(G_{\lambda})$ in the study of the charged plates in electrostatic
actuators. They showed that there exists $0<\lambda^{*}<\infty$ such
that for $\lambda \in(0,\lambda^{*})$ $(G_{\lambda})$ has a minimal
regular solutions $u_{\lambda}$ ($\sup_{B}u_{\lambda}<1$) while for
$\lambda>\lambda^{*}$, $(G_{\lambda})$ does not have any regular
solution. Moreover, the branch $\lambda\rightarrow u_{\lambda}(x)$
is increasing for each $x \in B$, and therefore the function
$u^{*}=\lim_{\lambda\nearrow \lambda^{*}}u_{\lambda}$ can be
considered as a generalized solution that corresponds to the pull-in
voltage $\lambda^{*}$. Now the important question is whether the
extremal solution $u^{*}$ is regular or not. In a recent paper Guo
and Wei \cite{GW} proved that the extremal solution $u^{*}$ is
regular for dimensions $N\leq 4$. In this paper we consider the
problem $(G_{\lambda})$ on the unit ball in $\R^N$:
$$ \left \{ \begin{array}{cl}
\beta \Delta^2 u- \tau \Delta u=\frac{\lambda}{(1-u)^2} \;\;
&\mbox{in $B$},\\ 0<u \leq 1\;\; &\mbox{in $B$},\\
u=\Delta u=0 \;\; &\mbox{on $\partial B$,}
\end{array} \right. \eqno(P_\lambda)$$
and show that the critical dimension for $(P_{\lambda})$ is $N=9$.
Indeed we prove that the extremal solution of $(P_{\lambda})$ is
regular ($\sup_{B}u^{*}<1$) for $N\leq 8$ and $\beta, \tau>0$ and it
is singular ($\sup_{B}u^{*}=1$) for $N\geq 9$, $\beta>0$, and
$\tau>0$ with $\frac{\tau}{\beta}$ small. Our proof of regularity of
the extremal solution in dimensions $5\leq N\leq 8$ is heavily
inspired by \cite{CEG} and \cite{DDGM}. On the other hand we shall
use certain improved Hardy-Rellich inequalities to prove that the
extremal solution is singular in dimensions $N\geq9$. Our improve
Hardy-Rellich inequalities follow from the recent result of
Ghoussoub-Moradifam \cite{GM} about Hardy and Hardy-Rellich
inequalities.

We now start by recalling some of the results from \cite{GW}
concerning $(P_{\lambda})$ that will be needed in the sequel. Define
\[ \lambda^*(B):= \sup \{ \lambda > 0: (P_{\lambda}) \; \mbox{ has a classical solution} \}.\]
We now introduce the following  notion of solution.

 \begin{definition} We say that $u$ is a {\it weak solution} of $(G_{\lambda})$,  if $0 \leq u \le 1$ a.e. in $\Omega$, $ \frac{1}{(1-u)^2} \in L^1(\Omega)$ and if
\[
\int_\Omega u (\beta \Delta^2 \phi- \tau \Delta \phi)\,dx = \lambda
\int_\Omega \frac{\phi}{(1-u)^2}dx, \qquad \forall \phi \in
W^{4,2}(\Omega)\cap H^1_0(\Omega),
\]
Say that $ u $ is a {\it weak super-solution} (resp. {\it weak
sub-solution}) of $(G_{\lambda})$, if the equality is replaced with
$ \ge $ (resp. $ \le $) for $ \phi \ge 0$.
\end{definition}

We now introduce the notion of stability. First, we equip the
function space  $\mathcal{H}:=H^2(\Omega) \cap H_0^1
(\Omega)=W^{2,2}(\Omega) \cap H_0^1 (\Omega)$
 with the norm
\[
\|\psi\|=\Big(\int_\Omega [\tau |\nabla \psi|^2+\beta |\Delta
\psi|^2] dx \Big)^{1/2}.
\]

\begin{definition} We say that a weak solution $u_\lambda$ of $(G_\lambda)$ is stable (respectively semi-stable) if the first eigenvalue $\mu_{1,
\lambda} (u_\lambda)$ of the problem
\begin{equation}
\label{e.1} -\tau \Delta h+\beta \Delta^2 h-\frac{2
\lambda}{(1-u_\lambda)^3} h=\mu h \;\;\mbox{in $\Omega$}, \;\;\;
h=\Delta h=0 \;\; \mbox{on $\partial \Omega$}
\end{equation}
is positive (resp., nonnegative).
\end{definition}

The operator $\beta \Delta^2 u- \tau \Delta u$ satisfies the
following maximum principle which will be frequently used in the
sequel.

\begin{lemma}(\cite{GW})\label{MAX}
Let $u\in L^1(\Omega)$.  Then $u\geq 0$ a.e. in $\Omega$, provided
one of the following conditions hold:
\begin{enumerate}

\item $u\in C^4(\overline{\Omega})$,  $\beta \Delta^2 u- \tau \Delta u\geq 0$ on $\Omega$, and  $u=\Delta u= 0$ on $\partial \Omega$.
\item $\int_{\Omega} u (\beta \Delta^2 \phi- \tau \Delta \phi)\,dx\geq 0$ for all $0\leq \phi \in W^{4,2}(\Omega)\cap H^1_0(\Omega)$.
\item $u\in W^{2,2}(\Omega)$, $u=0$, $\Delta u \leq 0$ on $\partial B$, and
$ \int_{\Omega} \big[\beta \Delta u \Delta  \phi+\tau \nabla u\nabla
\phi\big] dx  \geq 0$  for all $0\leq \phi \in W^{2,2}(\Omega)\cap
H^1_0(\Omega)$.
\end{enumerate}
Moreover, either $u\equiv 0$ or $u>0$ a.e. in $\Omega$.
\end{lemma}

\section{The pull-in voltage}

As in \cite{DDGM} and \cite{CEG}, we are led here to examine problem
$(P_\lambda)$ with non-homogeneous boundary conditions such as
$$ \left \{ \begin{array}{cl}
\beta \Delta^2 u- \tau \Delta u=\frac{\lambda}{(1-u)^2} \;\;
&\mbox{in $B$},\\ \alpha<u \leq 1\;\; &\mbox{in $B$},\\
u=\alpha, \; \Delta u=\gamma \;\; &\mbox{on $\partial B$,}
\end{array} \right. \eqno(P_\lambda, \alpha, \gamma)$$
where $\alpha, \gamma$ are given. Whenever we need to emphasis the
parameters $\beta$ and $\tau$ we will refer to problem $(P_{\lambda,
\alpha, \gamma})$ as $(P_{\lambda,\beta,\tau, \alpha, \gamma})$. In
this section and Section 3 we will obtain several results for the
following general form of $(P_\lambda, \alpha, \gamma)$
$$ \left \{ \begin{array}{cl}
\beta \Delta^2 u- \tau \Delta u=\frac{\lambda}{(1-u)^2} \;\;
&\mbox{in $\Omega$},\\ \alpha<u \leq 1\;\; &\mbox{in $\Omega$},\\
u=\alpha, \; \Delta u=\gamma \;\; &\mbox{on $\partial \Omega$,}
\end{array} \right. \eqno(G_\lambda, \alpha, \gamma)$$
which are analogous to the results obtained by Gui and Wei for
$(G_{\lambda})$ in \cite{GW}.

Let $\Phi$ denote the unique solution of
\begin{equation} \label{Phi} \left\{ \begin{array}{ll}
\beta \Delta^2 \Phi -\tau \Delta \Phi= 0 &\hbox{in } \Omega, \\
\Phi = \alpha\:,\:\: \Delta \Phi = \gamma&\hbox{on }\partial \Omega.
\end{array}\right.
\end{equation}
We will say that the pair $ (\alpha, \gamma)$ is admissible if
$\gamma \leq 0$, $\alpha < 1$, and $\sup_{\Omega}\Phi<1$. We now
introduce a notion of weak solution.
\begin{definition} We say that $u$ is a {\it weak solution} of $(P_{\lambda,\alpha,\gamma})$,  if $\alpha \leq u \le 1$ a.e. in $\Omega$, $ \frac{1}{(1-u)^2} \in L^1(\Omega)$ and if
\[
\int_\Omega (u-\Phi)(\beta\Delta^2 \phi-\tau\Delta \phi) = \lambda
\int_\Omega \frac{\phi}{(1-u)^2} \qquad \forall \phi \in
W^{4,2}(\Omega) \cap H_0^1(\Omega),
\]
where $\Phi$ is given in (\ref{Phi}). We say $ u $ is a {\it weak
super-solution} (resp. {\it weak sub-solution}) of
$(P_{\lambda,\alpha,\gamma})$, if the equality is replaced with $
\ge $ (resp. $ \le $) for $ \phi \ge 0$.
\end{definition}
\begin{definition} We say a weak solution $u$ of $(P_{\lambda,\alpha,\gamma})$ is regular (resp. singular) if $\|u\|_\infty<1$ (resp. $\|u\|_\infty=1$).
\end{definition}
  \noindent We now define
\[ \lambda^*(\alpha,\gamma):= \sup \Big\{ \lambda > 0: (P_{\lambda, \alpha, \gamma}) \; \mbox{ has a classical solution} \Big\}
\]
and
\[
\lambda_*(\alpha, \gamma):= \sup \Big\{ \lambda > 0: (P_{\lambda,
\alpha, \gamma}) \; \mbox{ has a weak solution} \Big\}.
\]
Observe that by the Implicit Function Theorem,  we can classically
solve $(P_{\lambda,\alpha,\gamma})$ for small $\lambda$'s.
Therefore, $\lambda^*(\alpha, \gamma)$ and  $\lambda_*(\alpha,
\gamma)$ are well defined for any admissible pair $(\alpha,
\gamma)$. To cut down on notations we won't always indicate $ \alpha
$ and $ \gamma$.  For example, $\lambda_* $ and $ \lambda^*$ will
denote the ``weak and strong critical voltages" of $ (P_{\lambda,
\alpha , \gamma})$.

Now let $U$ be a weak super-solution of
$(P_{\lambda,\alpha,\gamma})$ and recall the following existence
result.
\begin{theorem}(\cite{GW}) \label{exist} For every $0\leq f \in L^1(\Omega)$ there exists a unique $0\leq u \in L^1(\Omega)$ which satisfies
$$
\int_\Omega u (\beta \Delta^2 \phi- \tau \Delta
\phi)\,dx=\int_\Omega f \phi\, dx,
$$
for all $\phi \in W^{4,2}(\Omega)\cap H^1_0(\Omega)$.
\end{theorem}

We can introduce the following ``weak" iterative scheme: $u_0=U$ and
(inductively) let $u_n$, $n \geq 1$, be the solution of
$$\int_\Omega (u_n-\Phi) (\beta\Delta^2 \phi-\tau \Delta \phi)=\lambda \int_\Omega  \frac{\phi}{(1-u_{n-1})^2}\qquad\:\forall \: \phi \in W^{4,2}(\bar \Omega) \cap H_0^1(\Omega)$$
given by Theorem \ref{exist}. Since $0$ is a sub-solution of
$(P_{\lambda,\alpha,\gamma})$, inductively it is easily shown by
Lemma \ref{MAX} that $\alpha \leq u_{n+1}\leq u_n \leq U$ for every
$n \geq 0$. Since
$$(1-u_n)^{-2}\leq (1-U)^{-2} \in L^1(\Omega),$$
by Lebesgue Theorem the function $u=\displaystyle \lim_{n \to
+\infty} u_n$ is a weak solution of $(P_{\lambda,\alpha,\gamma})$ so
that $\alpha \leq u\leq U$. We therefore have the following result.

\begin{lemma} \label{super} Assume the existence of a weak super-solution $U$ of $(P_{\lambda,\alpha,\gamma})$.
 Then there exists a weak solution $u$ of $(P_{\lambda,\alpha,\gamma})$ so that $\alpha \leq u \leq U$ a.e. in $\Omega$.
\end{lemma}
\noindent In particular, for every $\lambda \in (0,\lambda_*)$, we
can find a weak solution of $(P_{\lambda,\alpha,\gamma})$. In the
same range of $\lambda'$s, this is still true for regular weak
solutions as shown in the following lemma.

\begin{lemma} \label{cch} Let $ (\alpha, \gamma)$ be an admissible pair and $u$ be a weak solution of $(P_{\lambda,\alpha,\gamma})$.
Then,  there exists a regular solution for every $0<\mu<\lambda$.
\end{lemma}
\noindent{\bf Proof:}  Let $\epsilon \in (0,1)$ be given and let  $
\bar u=(1-\epsilon)u+\epsilon \Phi$, where $\Phi$ is given in
(\ref{Phi}). By Lemma \ref{MAX} $\sup_{\Omega}
\Phi<\sup_{\Omega}u\leq 1$. Hence
$$
\sup_\Omega \bar u\leq (1-\epsilon)+\epsilon \sup_\Omega
\Phi<1\:,\quad \inf_\Omega \bar u\geq (1-\epsilon)\alpha +\epsilon
\inf_\Omega \Phi=\alpha,
$$
and for every $0\leq \phi \in W^{4,2}(\bar \Omega) \cap
H_0^1(\Omega)$ there holds:
\begin{eqnarray*}
\int_\Omega (\bar u-\Phi) (\beta\Delta^2 \phi-\tau \Delta \phi) &=&
(1-\epsilon) \int_\Omega (u-\Phi)(\beta\Delta^2 \phi-\tau \Delta
\phi)
= (1-\epsilon)\lambda \int_\Omega \frac{\phi}{(1-u)^2}\\
&=& (1-\epsilon)^3 \lambda \int_\Omega \frac{\phi}{(1-\bar
u+\epsilon (\Phi-1))^2} \geq (1-\epsilon)^3 \lambda \int_\Omega
\frac{\phi}{(1-\bar u)^2}.
\end{eqnarray*}
Note that $ 0 \le (1-\epsilon)(1-u)=1 - \bar{u}+\epsilon (\Phi -1)
<1-\bar u$. So $ \bar{u}$ is a weak super-solution of $ (P_{
(1-\epsilon)^3 \lambda, \alpha , \gamma})$ so that $\displaystyle
\sup_\Omega \bar u<1$. By  Lemma \ref{super} we get the existence of
a weak solution $w$ of $ (P_{ (1-\epsilon)^3 \lambda, \alpha ,
\gamma})$ so that $\alpha \leq w\leq \bar u$. In particular,
$\displaystyle \sup_\Omega w<1$ and $w$ is a regular weak solution.
Since $\epsilon \in (0,1)$ is arbitrarily chosen, the proof is done.
$\Box$
\endproof

Lemma \ref{cch} implies the existence of a regular weak solution
$U_\lambda$ for every $\lambda \in (0,\lambda_*)$. Introduce now a
``classical" iterative scheme: $u_0=0$ and (inductively)
$u_n=v_n+\Phi$, $n \geq 1$, where $v_n \in W^{4,2}(\Omega)\cap
H^1_0(\Omega)$ is the solution of
\begin{equation} \label{pranzo}
\beta \Delta^2 v_n- \tau \Delta  v_n=\beta \Delta^2 u_n- \tau \Delta
u_n= \frac{\lambda}{(1-u_{n-1})^2} \quad\hbox{in }\Omega \ \
\hbox{and} \ \ \Delta v_{n}=0 \quad\hbox{on } \partial \Omega.
\end{equation}
Since $v_n \in W^{4,2}(\Omega)\cap H^1_0(\Omega)$, $u_n$ is also a
weak solution of (\ref{pranzo}), and by Lemma \ref{MAX} we know that
$\alpha \leq u_n\leq u_{n+1} \leq U_\lambda$ for every $n \geq 0$.
Since $\displaystyle \sup_\Omega u_n \leq \displaystyle \sup_\Omega
U_\lambda<1$ for $n\geq 0$, we get that $(1-u_{n-1})^{-2} \in
L^2(\Omega)$ and the existence of $v_n$ is guaranteed. Since $v_n$
is easily seen to be uniformly bounded in $H^2(\Omega)$, we have
that $u_\lambda:=\displaystyle \lim_{n \to +\infty}u_n$ does hold
pointwise and weakly in $H^2(\Omega)$. By Lebesgue theorem, we have
that $u_\lambda$ is a radial weak solution of $(P_{\lambda})$ so
that $\displaystyle \sup_\Omega u_\lambda\leq \displaystyle
\sup_\Omega U_\lambda<1$. By elliptic regularity theory \cite{ADN},
$u_\lambda \in C^\infty(\bar \Omega)$ and $u_\lambda=\Delta
u_\lambda=0$ on $\partial \Omega$. So we can integrate by parts to
get
$$
\int_\Omega \beta (\Delta^2 u_\lambda- \tau \Delta  u_\lambda) \phi
\, dx =\int_\Omega u_\lambda (\beta \Delta^2 \phi- \tau \Delta
\phi)\, dx=\lambda \int_\Omega \frac{\phi}{(1-u_\lambda)^2}$$ for
every $\phi \in W^{4,2}(\Omega)\cap H^1_0(\Omega)$. Hence,
$u_\lambda$ is a classical solution of $(P_{\lambda})$ showing that
$\lambda^*=\lambda_*$.

Since the argument above shows that $u_\lambda<U$ for any other
classical solution $U$ of $(P_{\mu},\alpha,\gamma)$ with $\mu \geq
\lambda$, we have that $u_\lambda$ is exactly the minimal solution
and $u_\lambda$ is strictly increasing as $\lambda \uparrow
\lambda^*$. In particular, we can define $ u^*$ in the usual way: $
u^*(x)= \displaystyle \lim_{\lambda \nearrow \lambda^*}
u_\lambda(x)$.

\begin{lemma}  $\lambda^* (\Omega)<+\infty$.
\end{lemma}
\noindent{\bf Proof:}  Let $u$ be a classical solution of $
(P_{\lambda,\alpha, \gamma})$  and let $ (\psi, \mu_1)$ with $\Delta
\psi=0$ on $\partial \Omega$ denote the first eigenpair of $\beta
\Delta^2 - \tau \Delta $ in $H^2(\Omega)\cap H_0^1(\Omega)$ with $
\psi
>0$. Now let $C$ be such that
\[\int_{\partial \Omega}\left((\tau \alpha-\beta \gamma)\partial_{\nu}\psi-\beta \alpha \partial_{\nu}(\Delta \psi)\right)=C\int_\Omega \psi.\]

Multiplying $ (P_{\lambda,\alpha, \gamma})$ by $ \psi$ and then
integrating by parts one arrives at
\[ \int_\Omega \left( \frac{ \lambda}{(1-u)^2} - \mu_1 u -C\right) \psi =0. \]
Since $ \psi>0$ there must exist a point $\bar x \in \Omega$ where
$\frac{ \lambda}{(1-u(\bar x))^2} - \mu_1 u(\bar x)-C \le 0$. Since
$\alpha <u(\bar x)<1$, hence one can conclude that $ \lambda \le
\sup_{0< u <1} ( \mu_1 u+C)(1-u)^2$, which shows that $
\lambda^*<+\infty$. $\Box$\\
\endproof

In conclusion, we have shown the following description of the
minimal branch.

\begin{theorem}  $\lambda^* \in (0,+\infty)$ and the following holds:
\begin{enumerate}
\item For each $ 0 < \lambda < \lambda^*$ there exists a regular and minimal solution $u_\lambda$ of $(P_{\lambda,\alpha,\gamma})$.
\item For each $ x \in \Omega$ the map $ \lambda \mapsto u_\lambda(x)$ is strictly increasing on $ (0,\lambda^*)$.
\item For $ \lambda > \lambda^*$ there are no weak solutions of $(P_{\lambda,\alpha, \gamma})$.
\end{enumerate}
\label{quasi}
\end{theorem}

\section{Stability of the minimal solutions}
This section is devoted to the proof of the following stability
result for minimal solutions. We shall need  the following notion of
$\mathcal{H}-$weak solutions, which is an intermediate class between
classical and weak solutions.

\begin{definition} We say that $u$ is an $\mathcal{H}-$weak solution of $(P_{\lambda,\alpha, \gamma})$ if $u-\Phi \in H^2(\Omega)\cap H_0^1(\Omega)$, $0 \le u \le 1$ a.e. in $\Omega$, $ \frac{1}{(1-u)^2} \in L^1(\Omega)$ and
\[
 \int_\Omega \big[\beta \Delta u \Delta  \phi+\tau \nabla u\nabla   \phi\big] dx
= \lambda \int_\Omega \frac{\phi}{(1-u)^2}, \qquad \forall \phi \in
W^{2,2}(\Omega)\cap H^1_0(\Omega).
\]
where $\Phi$ is given by (\ref{Phi}). We say that $u$ is an
$\mathcal{H}-$weak super-solution (resp. an $\mathcal{H}-$weak
sub-solution) of $(P_{\lambda, \alpha, \gamma})$ if for $ \phi \ge
0$ the equality is replaced with $ \ge$ (resp. $ \le $) and $u\geq
0$ (resp. $\leq$), $\Delta u \leq 0$ (resp. $\geq$) on $\partial
\Omega$.
\end{definition}

\begin{theorem} \label{stable.square}  Suppose that $(\alpha, \gamma)$
is an admissible pair.
\begin{enumerate}
\item  The minimal solution $ u_\lambda $ is stable, and is the unique semi-stable $\mathcal{H}-$weak solution of $(P_{\lambda,\alpha,\gamma})$.
\item The function $ u^*:= \displaystyle \lim_{\lambda \nearrow \lambda^*} u_\lambda$ is a well-defined semi-stable $\mathcal{H}-$weak solution of  $(P_{\lambda^*,\alpha,\gamma})$.
\item $u^*$ is the unique $\mathcal{H}-$weak solution of $(P_{\lambda^*,\alpha,\gamma})$, and when $u^*$ is classical solution, then $\mu_1(u^*)=0$.
\item If $ v$ is a singular, semi-stable $\mathcal{H}-$weak solution of $ (P_{ \lambda,\alpha,\gamma})$, then $ v=u^*$ and $ \lambda = \lambda^*$
\end{enumerate}
\end{theorem}
\noindent The main tool is the following  comparison lemma which is
valid exactly in the class $\mathcal{H}$.
\begin{lemma} \label{shi} Let $ (\alpha, \gamma)$ be an admissible pair and $u$ be a semi-stable $\mathcal{H}-$weak solution of $(P_{\lambda, \alpha, \gamma})$.
Assume $ U $ is a $\mathcal{H}-$weak super-solution of $(P_{\lambda,
\alpha, \gamma})$. Then
\begin{enumerate}
\item $ u \le U$ a.e. in $\Omega$;
\item If $ u$ is a classical solution and $ \mu_1(u)=0$ then $ U=u$.
\end{enumerate}
\end{lemma}
\noindent{\bf Proof:}  (i)   Define $ w:= u-U$.   Then by means of
the Moreau decomposition for the biharmonic operator (see \cite{M}
and \cite{B}), there exist $ w_1$ and $w_2 \in H^2(\Omega)\cap
H^{1}_{0}(\Omega)$, with $ w=w_1 + w_2$, $ w_1 \ge 0$ a.e., $\beta
\Delta^2 w_2-\tau \Delta w_{2} \le 0 $ in the $\mathcal{H}-$weak
sense and
$\int_\Omega \beta \Delta w_1 \Delta w_2+\tau \nabla w_1.\nabla w_2=0$.   Lemma \ref{MAX} gives that $w_2 \le 0$ a.e. in $\Omega$.\\
Given $ 0 \le \phi \in C_c^\infty(\Omega)$,  we have
\[ \int_\Omega \beta \Delta w \Delta \phi+\tau\nabla w.\nabla \phi \leq  \lambda \int_\Omega (f(u) - f(U)) \phi, \]
where $ f(u):= (1-u)^{-2}$.  Since $ u$ is semi-stable,  one has
\begin{eqnarray*}
\lambda \int_\Omega f'(u) w_1^2 \le  \int_\Omega \beta(\Delta
w_1)^2+\tau |\nabla w_{1}|^2 = \int_\Omega \beta \Delta w \Delta
w_1+\tau \nabla w.\nabla w_1 \le \lambda \int_\Omega ( f(u) - f(U))
w_1.
\end{eqnarray*}
Since $ w_1 \ge w$ one has
\[ \int_\Omega f'(u) w w_1 \le \int_\Omega (f(u)-f(U)) w_1,\]
which re-arranged gives
\[ \int_\Omega \tilde{f} w_1 \geq 0,\]
where $ \tilde{f}(u)= f(u) - f(U) -f'(u)(u-U)$. The strict convexity
of $f$ gives $ \tilde{f} \le 0$ and $ \tilde{f}< 0 $ whenever $u
\not= U$. Since $w_1 \ge 0$ a.e. in $\Omega$, one sees that $ w \le
0 $ a.e.  in $\Omega$. The inequality $ u \le U$  a.e. in $\Omega$
is then established.

\medskip \noindent (ii) Since $u$ is a classical solution, it is easy to see that the infimum of $\mu_1(u)$ is attained at some $\phi$.
The function $\phi$ is then the first eigenfunction of $\beta
\Delta^2-\tau \Delta-\frac{2\lambda}{(1-u)^3}$ in $H^2(\Omega)\cap
H^1_0(\Omega)$. Now we show that $ \phi$ is of fixed sign.  Using
the above decomposition, one has $ \phi= \phi_1 + \phi_2$ where $
\phi_i \in H^2(\Omega)\cap H^1_0(\Omega)$ for $i=1,2$, $ \phi_1 \ge
0$, $ \int_\Omega \beta \Delta \phi_1 \Delta \phi_2+\tau \nabla
\phi_{1}.\nabla \phi_{2}=0$ and $\beta \Delta^2 \phi_2 -\tau \Delta
\phi_2\le 0$ in the $\mathcal{H}-$weak sense. If $ \phi$ changes
sign, then $ \phi_1 \not\equiv 0$ and $ \phi_2 <0$ in $\Omega$
(recall that either $\phi_2<0$ or $\phi_2=0$ a.e. in $\Omega$). We
can write now
\begin{eqnarray*}
0 = \mu_1(u) &\le&  \frac{ \int_\Omega \beta( \Delta (\phi_1
-\phi_2))^2+\tau |\nabla(\phi_1 -\phi_2)|^2 - \lambda f'(u) ( \phi_1
- \phi_2)^2}{ \int_\Omega ( \phi_1 - \phi_2)^2}\\ &<& \frac{
\int_\Omega \beta( \Delta \phi)^2+\tau |\nabla \phi|^2 - \lambda
f'(u) \phi^2 }{ \int_\Omega \phi^2} =\mu_1(u),
\end{eqnarray*}
in view of $\phi_1 \phi_2<-\phi_1\phi_2$ in a set of positive
measure, leading to a contradiction.

So we can assume $ \phi \ge 0$, and by Lemma \ref{MAX} we have
$\phi>0$ in $\Omega$. For $ 0 \le t \le 1$,  define
$$g(t)=\int_\Omega \beta \Delta \left[t U+(1-t)u \right] \Delta
\phi+\tau\nabla \left[t U+(1-t)u \right].\nabla \phi  - \lambda
\int_\Omega f( tU+(1-t)u) \phi,$$ where $\phi$ is the above first
eigenfunction. Since $ f$ is convex one sees that
$$g(t)\geq \lambda \int_\Omega \left[t f(U)+(1-t)f(u)-f(tU+(1-t)u)\right]\phi \geq 0$$  for every $t \geq 0$. Since $ g(0) =0$ and
$$ g'(0)= \int_\Omega \beta \Delta (U-u) \Delta \phi +\tau \nabla(U-u).\nabla \phi-\lambda f'(u)(U-u)\phi=0 ,$$
we get that
\[ g''(0)=- \lambda \int_\Omega f''(u) (U-u)^2 \phi\geq 0.\]
Since $f''(u)\phi>0$ in $\Omega$, we finally get that $ U=u$ a.e. in
$\Omega$.
$\Box$\\
\endproof

A more general version of Lemma \ref{shi} is available in the
following.
\begin{lemma} \label{poo} Let $ (\alpha,\gamma)$ be an admissible pair and $\gamma'\leq 0$. Let $u$ be a semi-stable $\mathcal{H}-$weak sub-solution of
$(P_{\lambda, \alpha,\gamma})$ with $u=\alpha'\leq \alpha$, $\Delta
u=\beta' \geq \beta$ on $\partial \Omega$. Assume that $U$ is a
$\mathcal{H}-$weak super-solution of $(P_{\lambda, \alpha ,\gamma})$
with $U=\alpha$, $\Delta U=\beta$ on $\partial \Omega$. Then $ U \ge
u$ a.e. in $\Omega$.
\end{lemma}
\noindent{\bf Proof:}   Let $ \tilde{u} \in H^2(\Omega)\cap
H^1_0(\Omega)$ denote a weak solution of $ \beta \Delta^2
\tilde{u}-\tau \Delta\tilde{u}= \beta\Delta^2 (u-U)-\tau
\Delta(u-U)$ in $\Omega$ and $\tilde{u}=\Delta \tilde{u}=0$ on
$\partial \Omega$. Since $\tilde{u}-u+U\geq0$ and
$\Delta(\tilde{u}-u+U)\leq 0$ on $\partial \Omega$, by Lemma
\ref{MAX} one has that $ \tilde{u} \ge u-U $ a.e. in $\Omega$. By
means of the Moreau decomposition (see \cite{M} and \cite{B}) we
write $\tilde u$ as $ \tilde{u} = w+v $, where $ w,v \in
H_0^2(\Omega)$, $ w \ge 0 $ a.e. in $\Omega$, $ \beta\Delta^2 v-\tau
\Delta v \le 0$ in a $\mathcal{H}-$weak sense and $\int_\Omega \beta
\Delta w \Delta v+\tau \nabla w.\nabla v=0$. Then for $ 0 \le \phi
\in W^{4,2} (\bar \Omega)\cap H_0^1(\Omega)$, one has
\[ \int_\Omega \beta \Delta \tilde{u} \Delta \phi +\tau \nabla \tilde{u}. \nabla \phi \leq \lambda \int_\Omega (f(u)- f( U)) \phi .\]
In particular, we have
\[ \int_\Omega \beta \Delta \tilde{u} \Delta w +\tau \nabla \tilde{u}. \nabla w\le \lambda \int_\Omega ( f(u)-f(U)) w.\]
Since the semi-stability of $u$ gives that
\begin{eqnarray*}
\lambda \int_\Omega f'(u) w^2\leq \int_\Omega \beta( \Delta
w)^2+\tau|\nabla w|^2 = \int_\Omega \beta\Delta \tilde{u} \Delta
w+\tau \nabla \tilde{u}.\nabla w ,
\end{eqnarray*}
we get that
\[  \int_\Omega f'(u) w^2 \le \int_\Omega ( f(u)-f(U)) w.\]  By Lemma \ref{MAX} we have $v\leq 0$ and then $ w \ge \tilde{u} \ge u -U$ a.e. in $\Omega$. So we obtain that
\[ 0 \le \int_\Omega \left( f(u)-f(U)-f'(u)(u-U) \right) w.\]  The strict convexity of $ f$ implies that $ U \ge u $ a.e. in
$\Omega$. $\Box$ \\
\endproof

We need also some a-priori estimates along the minimal branch
$u_\lambda$.

\begin{lemma} \label{extremalsol.10} Let $(\alpha,\gamma)$ be an admissible pair. Then for every   $\lambda \in (0,\lambda^*)$, we have
 \[
 2 \int_\Omega \frac{(u_\lambda-\Phi)^2}{(1-u_\lambda)^3} \le \int_\Omega \frac{ u_\lambda-\Phi}{(1-u_\lambda)^2},
 \]
where $\Phi$ is given by (\ref{Phi}). In particular, there is a
constant $C>0$ independent of $\lambda$ so that
\begin{equation} \label{tardi.10}
\int_\Omega (\tau |\nabla u_\lambda |^2+\beta |\Delta u_\lambda
|^2)dx+\int_\Omega \frac{1}{(1-u_\lambda)^3} \leq C,
\end{equation}
for every $\lambda \in (0,\lambda^*)$.
\end{lemma}

\noindent{\bf Proof:}  Testing $ (P_{\lambda,\alpha,\gamma})$ on $
u_\lambda-\Phi \in W^{4,2}( \Omega) \cap H^1_0(\Omega)$, we see that
\begin{eqnarray*}
\lambda \int_\Omega \frac{ u_\lambda-\Phi}{(1-u_\lambda)^2} =
\int_\Omega (\tau |\nabla ( u_\lambda-\Phi) |^2+\beta (\Delta (
u_\lambda-\Phi) )^2)dx \ge 2 \lambda \int_\Omega \frac{ (
u_\lambda-\Phi)^2}{( 1-u_\lambda)^3}.
\end{eqnarray*}
In the view of $\beta\Delta^2 \Phi-\tau \Delta \Phi=0$. In
particular, for $\delta>0$ small we have that
\begin{eqnarray*}
\int_{\{|u_\lambda| \geq \delta \}}\frac{1}{(1-u_\lambda)^3}&\leq &
\frac{1}{\delta^2} \int_{\{|u_\lambda-\Phi| \geq \delta \}}\frac{(
u_\lambda-\Phi)^2}{(1-u_\lambda)^3}
\leq \frac{1}{\delta^2} \int_\Omega \frac{1}{(1-u_\lambda)^2}\\
&\leq &\delta \int_{\{|u_\lambda-\Phi| \geq \delta
\}}\frac{1}{(1-u_\lambda)^3}+C_\delta
\end{eqnarray*}
by means of Young's inequality. Since for $\delta$ small
$$\int_{\{|u_\lambda-\Phi| \leq \delta \}}\frac{1}{(1-u_\lambda)^3}\leq C,
$$
for some $C>0$, we get that
$$
\int_\Omega \frac{1}{(1-u_\lambda)^3} \leq C,
$$
for some $C>0$ and for every $\lambda \in (0,\lambda^*)$. Since
\begin{eqnarray*}
\int_\Omega (\tau |\nabla u_\lambda |^2+\beta |\Delta u_\lambda |^2)dx&=& \int_{\Omega}\left (\beta \Delta u_{\lambda}\Delta \Phi+\tau\nabla u_{\lambda}.\nabla\Phi\right) +\lambda \int_\Omega \frac{u_\lambda-\Phi}{(1-u_\lambda)^2}\\
&\leq& \delta \int_\Omega (\tau |\nabla u_\lambda |^2+\beta |\Delta
u_\lambda |^2)dx +C_\delta+C \left(\int_\Omega
\frac{1}{(1-u_\lambda)^3} \right)^{\frac{2}{3}}
\end{eqnarray*}
in view of Young's and H\"older's inequalities, estimate (\ref{tardi.10}) is finally established. $\Box$\\

\noindent {\bf Proof of Theorem \ref{stable.square}:} (1)\, Since
$\|u_\lambda\|_\infty <1$, the infimum defining $\mu_1(u_\lambda)$
is achieved at a first eigenfunction for every $\lambda \in
(0,\lambda^*)$. Since $\lambda  \mapsto u_\lambda(x)$ is increasing
for every $x \in \Omega$, it is easily seen that $\lambda \mapsto
\mu_1( u_\lambda)$ is a decreasing and  continuous function on $ (0,
\lambda^*)$.  Define
\[ \lambda_{**}:= \sup\{ 0 <\lambda < \lambda^*: \: \mu_1( u_\lambda) >0 \} .\]
We have that $ \lambda_{**}= \lambda^*$. Indeed, otherwise we would
have $ \mu_1(u_{ \lambda_{**}}) =0$, and for every $ \mu \in (
\lambda_{**}, \lambda^*)$,  $ u_{\mu}$ would be a classical
super-solution of $ (P_{ \lambda_{**},\alpha, \gamma})$. A
contradiction arises since Lemma \ref{shi} implies $u_{\mu} =
u_{\lambda_{**}}$.
Finally, Lemma \ref{shi} guarantees the uniqueness in the class of semi-stable $\mathcal{H}-$weak solutions.\\

(2) \, It follows from  (\ref{tardi.10}) that $u_\lambda \to u^*$ in
a pointwise sense and weakly in $H^2(\Omega)$, and $ \frac{1}{1-u^*}
\in L^3(\Omega)$. In particular, $u^*$ is a $H^2-$weak solution of
$(P_{ \lambda^*, \alpha ,\gamma})$ which is also semi-stable as the
limiting function of the
semi-stable solutions $\{u_\lambda\}$.\\

(3) Whenever $\|u^*\|_\infty<1$, the function $u^*$ is a classical solution, and by the Implicit Function Theorem  we have that $\mu_1(u^*)=0$ to
 prevent the continuation of the minimal branch beyond $\lambda^*$. By Lemma \ref{shi},  $u^*$ is then the unique $\mathcal{H}-$weak solution of
 $(P_{\lambda^*,\alpha,\gamma})$. \\

(4) \, If $ \lambda < \lambda^*$, we get by uniqueness that
$v=u_\lambda $. So $v$ is not singular and a contradiction arises.
Now, by  Theorem \ref{quasi}(3) we have that $ \lambda = \lambda^*$.
Since $ v $ is a semi-stable $\mathcal{H}-$ weak solution of $ (P_{
\lambda^*, \alpha, \gamma})$ and $ u^*$ is a $\mathcal{H}-$ weak
super-solution of $ (P_{\lambda^*, \alpha , \gamma})$, we can apply
Lemma \ref{shi} to get $ v \le u^*$ a.e. in $\Omega$.  Since $u^*$
is also a semi-stable solution, we can reverse the roles of $ v$ and
$ u^*$ in Lemma \ref{shi} to see that $ v \ge u^*$ a.e. in $\Omega$.
So equality $v=u^*$ holds and the proof is done. $\Box$

\section{Regularity of the extremal solutions in dimensions $N\leq8$}
In this section we shall show that the extremal solution is regular
in small dimensions. Let us begin with the following lemma.
\begin{lemma}\label{bclemm}
Let $N\geq 5$ and $(u^{*}, \lambda^{*})$ be the extremal pair of
$(P_{\lambda})$. If $u^{*}$ is singular, and he set
\begin{equation}
\Gamma:=\{r\in(0,1): u_{\delta}(r)>u^{*}(r)\}
\end{equation}
is non-empty, where $u_{\delta}(x):=1-C_{\delta}|x|^{\frac{4}{3}}$
and $C_{\delta}>1$ is a constant. Then there exists $r_{1} \in
(0,1)$ such that $u_{\delta}(r_{1})\geq u^{*}(r_{1})$ and $\Delta
u_{\delta}(r_{1})\leq \Delta u^{*}(r_{1})$.
\end{lemma}
{\bf Proof.} Assume by contradiction that for every $r$ with
$u_{\delta}(r_{1})\geq u^{*}(r_{1})$ one has $\Delta
u_{\delta}(r_{1})> \Delta u^{*}(r_{1})$. Since $\Gamma$ is non-empty
and
\[u_{\delta}(1)=1-C_{\delta}<0=u^{*}(1),\]
there exists $s_{1} \in (0,1)$ such that
$u_{\delta}(s_{1})=u^{*}(s_{1})$. We claim that
\[u_{\delta}(s)> u^{*}(s),\]
for $0<s< s_{1}$. Assume that there exist $s_{3}<s_{2}\leq s_{1}$
such that $u^{*}(s_{2})=u_{\delta}(s_{2})$,
$u^{*}(s_{3})=u_{\delta}(s_{3})$ and $u_{\delta}(s)\geq u^{*}(s)$
for $s \in (s_{3},s_{2})$. By our assumption $\Delta u_{s}> \Delta
u^{*}(s)$ for $s \in (s_{3},s_{2})$ which contradicts the maximum
principle and justifies the claim. Therefore
$u_{\delta}(s)>u^{*}(s)$ for $0<s<s_{1}$.  Now set
$w:=u_{\delta}-u^{*}$. Then $w\geq0$ on $B_{s_{1}}$ and $\Delta
w\leq 0$ in $B_{s_{1}}$. Since $w(0)=0$, by strong maximum principle
we get $w\equiv 0$ on $B_{s_{1}}$. This is a contradiction and
completes the proof. $\Box$

\begin{theorem} \label{bestim}
Let $N\geq 5$ and $(u^{*}, \lambda^{*})$ be the extremal pair of
$(P_{\lambda})$. When $u^{*}$ is singular, then
\begin{equation*}
1-u^{*}\leq C|x|^{\frac{4}{3}} \ \ in \ \ B,
\end{equation*}
where $C:=(\frac{\lambda^{*}}{\beta \bar{\lambda}})^{\frac{1}{3}}$
and $\bar{\lambda}:=\frac{8(N-\frac{2}{3})(N-\frac{8}{3})}{9}$.
\end{theorem}
{\bf Proof.} For $\delta>0$, define
$u_{\delta}(x):=1-C_{\delta}|x|^{\frac{4}{3}}$ with
$C_{\delta}:=(\frac{\lambda^{*}}{\beta
\bar{\lambda}}+\delta)^{\frac{1}{3}}>1$. Since $N\geq5$, we have
that $u_{\delta} \in H^{2}_{loc}(\R^{N})$ and $u_{\delta}$ is a
$\mathcal{H}-$weak solution of

\begin{equation*}
\beta \Delta^2 u_{\delta}-\tau \Delta
u_{\delta}=\frac{\lambda^{*}+\beta \delta
\bar{\lambda}}{(1-u_{\delta})^2}+\frac{4}{3}\tau
C_{\delta}(N-\frac{2}{3})|x|^{-\frac{2}{3}} \ \ in \ \ \R^N.
\end{equation*}
We claim that $u_{\delta}\leq u^{*}$ in $B$, which will finish the
proof by just letting $\delta \rightarrow 0$.

Assume by contradiction that the set $\Gamma:=\{r\in(0,1):
u_{\delta}(r)>u^{*}(r)\}$ is non-empty. By Lemma \ref{bclemm} the
set \begin{equation*} \Lambda:= \{r \in (0,1): u_{\delta}(r)\geq
u^{*}(r) \ \ \hbox{and} \ \ \Delta u_{\delta}(r)\leq \Delta
u^{*}(r)\}
\end{equation*}
is non-empty.  Let $r_{1} \in \Lambda$. Since
\begin{equation*}
u_{\delta}(1)=1-C_{\delta}<0=u^{*}(1),
\end{equation*}
we have that $0<r_{1}<1$. Define
\[\alpha:=u_{*}(r_{1})\leq u_{\delta}(r_{1}), \ \ \gamma := \Delta u^{*}(r_{1})\geq \Delta u_{\delta}(r_{1}).\]
Setting
$u_{\delta,r_{1}}=r_{1}^{-\frac{4}{3}}(u_{\delta}(r_{1}r)-1)+1$, we
see that $u_{\delta,r_{1}}$ is a $\mathcal{H}-$weak super-solution
of $(P_{\lambda^{*}+\delta \lambda,\beta,r_{1}^{-2}\tau,
\alpha',\gamma'})$, where
\begin{equation*}
\alpha':=r_{1}^{-\frac{4}{3}}(\alpha-1)+1, \ \
\gamma'=r_{1}^{\frac{2}{3}}\gamma.
\end{equation*}
Similarly, define
$u^{*}_{r_{1}}(r)=r_{1}^{-\frac{4}{3}}(u^{*}(r_{1}r)-1)+1$. Note
that $\Delta^2 u^{*}-\alpha \Delta u^{*} \geq 0$ in $B$ and $\Delta
u^{*}=0$ on $\partial B$. Hence by maximum principle we have $\Delta
u^{*}\leq 0$ in $B$ and therefore $\gamma'\leq 0$. Also obviously
$\alpha'<1$. So, $(\alpha', \gamma')$ is an admissible pair and by
Theorem \ref{stable.square}(4) we get that
$(u^{*}_{r_{1}},\lambda^{*})$ coincides with the extremal pair of
$(P_{ \lambda,\beta,r_{1}^{-2}\tau, \alpha',\gamma'})$ in $B$. Also
by Lemma \ref{super} we get the existence of a week solution of
$(P_{\lambda^{*}+\delta \lambda,\beta,r_{1}^{-2}\tau,
\alpha',\gamma'})$. Since $\lambda^{*}+\delta \lambda>\lambda^{*}$,
we contradict the fact that $\lambda^{*}$ is the extremal parameter
of $(P_{\lambda,\beta,r_{1}^{-2}\tau,
\alpha',\gamma'})$. $\Box$\\

Now we are ready to prove the following result.

\begin{theorem} If $5\leq N\leq 8$, then the extremal solution $u^*$
of $(P)_\lambda$ is regular.
\end{theorem}

{\bf Proof.} Assume that $u^{*}$ is singular. For $\epsilon>0$
define $\varphi(x):=|x|^{\frac{4-N}{2}+\epsilon}$ and note that
\begin{equation*}
(\Delta \varphi)^2=(H_{N}+O(\epsilon))|x|^{-N+2\epsilon}, \ \
\hbox{where} \ \ H_{N}:=\frac{N^2(N-4)^2}{16}.
\end{equation*}
Given $\eta \in C^{\infty}_{0}(B)$, and since $N\geq 5$, we can use
the test function $\eta \varphi \in H^{2}_{0}(B)$ into the stability
inequality to obtain
\begin{equation*}
2\lambda^{*}\int_{B}\frac{\varphi^2}{(1-u^{*})^{3}}\leq
\beta\int_{B}(\Delta \varphi)^2+\tau \int_{B} |\nabla \varphi|^2+
O(1),
\end{equation*}
where $O(1)$ is a bounded function as $\epsilon \rightarrow 0$. By
Theorem \ref{bestim} we find
\begin{equation*}
2\bar{\lambda}\int_{B}\frac{\varphi^2}{|x|^{4}}\leq \int_{B} (\Delta
\varphi)^2+O(1),
\end{equation*}
and then
\begin{equation*}
2\bar{\lambda}\int_{B}|x|^{-N+2\epsilon}\leq
(H_{N}+O(\epsilon))\int_{B}|x|^{-N+2\epsilon}+O(1).
\end{equation*}
Computing the integrals on obtains
\begin{equation*}
2\bar{\lambda}\leq H_{N}+O(\epsilon).
\end{equation*}
Letting $\epsilon\rightarrow 0$ we get $2\bar{\lambda}\leq H_{N}$.
Graphing this relation we see that $N\geq 9$. $\Box$

\section{The extremal solution is singular in dimensions $N\geq 9$}
In this section we will show that the extremal solution $u^{*}$ of
$(P_{\lambda, \beta, \tau, 0, 0})$ in dimensions $N\geq 9$ is
singular for $\tau>0$ sufficiently small. To do this, first we shall
show that the extremal solution of $(P_{\lambda, 1, 0, 0, 0})$ is
singular in dimensions $N\geq 9$. Again to cut down the notation we
won't always indicate that $\beta=1$ and $\tau=0$.

We have to distinguish between three different ranges for the
dimension. For each range, we will need a suitable Hardy-Rellich
type inequality that will be established in the appendix, by using
the recent results of Ghoussoub-Moradifam \cite{GM}.
\begin{trivlist}

\item {\bf $\bullet$ Case  $N\geq 16$:}\, To establish the singularity of $u^*$ for these dimensions we shall need the classical Hardy-Rellich
inequality, which is valid for all $\phi \in H^2(B)\cap H^1_0(B)$:
\begin{equation}\label{HR1}
 \int_B (\Delta \phi)^2\, dx \ge \frac{ N^2(N-4)^2}{16} \int_B \frac{\phi^2}{|x|^4} \, dx .
 \end{equation}

\item {\bf $\bullet$ Case  $10\leq N\leq 16$:}\,  For this case,  we shall need the following inequality valid  for all $\phi \in H^2(B)\cap H^1_0(B)$
\begin{eqnarray} \label{HR2}
\int_{B}(\Delta \phi)^2 &\geq&
\frac{(N-2)^2(N-4)^2}{16}\int_{B}\frac{\phi^2}{(|x|^2-\frac{N}{2(N-1)}|x|^{\frac{N}{2}+1})(|x|^2-|x|^{\frac{N}{2}})}\\
&&+\frac{(N-1)(N-4)^2}{4}
\int_{B}\frac{\phi^2}{|x|^2(|x|^2-|x|^{\frac{N}{2}})}.\nonumber
\end{eqnarray}

\item {\bf $\bullet$ Case  $N=9$:}\,  This case is the trickiest and  will require the following inequality  for all $\phi \in H^2(B)\cap H^1_0(B)$, which is valid for $N\geq 7$

\begin{equation}
\int_{B}|\Delta u|^2\geq \int_{B}W(|x|)u^2.
\end{equation}
where where
\begin{equation*}
W(r)=K(r)(\frac{(N-2)^2}{4(r^2-\frac{N}{2(N-1)}r^{\frac{N}{2}+1})}+
\frac{(N-1)}{r^2}),
\end{equation*}
\[K(r)=-\frac{\varphi''(r)+\frac{(n-3)}{r}\varphi'(r)}{\varphi(r)},\]
and
\[\varphi(r)=r^{-\frac{N}{2}+2}+9r^{-2}+10r-20.\]

\end{trivlist}

The next lemma will be our main tool to guarantee that $u^{*}$ is
singular for $N\geq 9$. The proof is based on an upper estimate by a
singular stable sub-solution.

\begin{lemma}\label{sing-lem}
Suppose there exist $\lambda'>0$ and a radial function $u \in
H^{2}(B)\cap W^{4,\infty}_{loc}(B\setminus \{0\})$ such that
\begin{equation}\label{cond1}
\Delta^2 u\leq \frac{\lambda'}{(1-u)^2} \ \ \hbox{for} \ \ 0<r<1,
\end{equation}

\begin{equation}
u(1)=0, \ \ \ \ \Delta u|_{r=1}=0,
\end{equation}
\begin{equation}
\hbox{u is singular},
\end{equation}
and
\begin{equation}\label{cond2}
2\beta \int_{B}\frac{ \varphi^2}{(1-u)^3}\leq \int_{B}(\Delta
\varphi)^2\ \ \hbox{for all} \ \ \varphi \in H^2(B)\cap H^1_0(B),
\end{equation}

for some $\beta>\lambda'$. Then $u^{*}$ is singular and
\begin{equation}\label{est1}
\lambda^{*} \leq \lambda'
\end{equation}
\end{lemma}
{\bf Proof.} By Lemma \ref{poo} we have (\ref{est1}). Let
$\frac{\lambda'}{\beta}<\gamma<1$ and
\begin{equation}
\alpha:=(\frac{\gamma \lambda^{*}}{\lambda'})^{1/3},
\end{equation}
and define $\bar{u}:=1-\alpha(1-u)$. We claim that
\begin{equation}\label{claim}
u^{*}\leq \bar{u} \ \ \hbox{in} \ \ B.
\end{equation}
To prove this, we shall show that for $\lambda<\lambda^{*}$
\begin{equation}
u_{\lambda}\leq \bar{u} \ \ \hbox{in} \ \ B.
\end{equation}
Indeed, we have
\begin{eqnarray*}
\Delta^2(\bar{u})=\alpha \Delta^2(\bar{u})\leq \frac{
\alpha\lambda'}{(1-u)^2}=\frac{\alpha^3\lambda'}{(1-\bar{u})^2}.
\end{eqnarray*}
By (\ref{est1}) and the choice of $\alpha$
\begin{equation*}
\alpha^3\lambda'<\lambda^{*}.
\end{equation*}
To prove (\ref{claim}) it suffices to prove it for
$\alpha^3\lambda'<\lambda<\lambda^{*}$. Fix such $\lambda$ and
assume that $(\ref{claim})$ is not true. Then
\begin{equation*}
\Lambda=\{0\leq R\leq 1\mid u_{\lambda}(R)>\bar{u}(R)\},
\end{equation*}
in non-empty. There exists $0<R_{1}<1$, such that
$u_{\lambda}(R_{1})\geq u^{*}(R_{1})$ and $\Delta
u_{\lambda}(R_{1})\leq \Delta u^{*}(R_{1})$, since otherwise we can
find $0<s_{1}<s_{2}<1$ so that $u_{\lambda}(s_{1})=\bar{u}(s_{1})$,
$u_{\lambda}(s_{2})=\bar{u}(s_{2})$, $ u_{\lambda}(R)>\bar{u}(R)$,
and $\Delta u_{\lambda}(R_{1})> \Delta u^{*}(R_{1})$ which
contradict the maximum principle. Now consider the following problem
\begin{eqnarray*}
\Delta^2 u&=&\frac{\lambda}{(1-u)^2} \ \ \hbox{in} \ \ B \\
u&=&u_{\lambda}(R_{1}) \ \ \hbox{on} \ \ \partial B \\
\Delta u&=&\Delta u_{\lambda} \ \ \hbox{on} \ \ \partial B.
\end{eqnarray*}
Then $u_{\lambda}$ is a solution to the above problem while
$\bar{u}$ is a sub-solution to the same problem. Moreover $\bar{u}$
is stable since,
\[\lambda<\lambda^{*}\]
and hence
\[\frac{2\lambda}{(1-\bar{u})^3}\leq \frac{2\lambda^{*}}{\alpha^3(1-u)^3}=\frac{2\lambda'}{ \gamma(1-u)^3}<\frac{2\beta}{(1-u)^3}.\]
We deduce $\bar{u}\leq u_{\lambda}$ in $B_{R_{1}}$ which is
impossible, since $\bar{u}$ is singular while $u_{\lambda}$ is
smooth. This establishes (\ref{claim}). From (\ref{claim}) and the
above two inequalities we have
\[\frac{2\lambda^{*}}{(1-u^{*})^3}\leq \frac{2\lambda'}{ \gamma(1-u)^3}<\frac{\beta}{(1-u)^3}.\]
Thus
\[\inf_{\varphi \in C^{\infty}_{0}}(B)\frac{\int_{B}(\Delta \varphi)^2-\frac{2\lambda^{*}\varphi^2}{(1-u^{*})^3}}{\int_{B}\varphi^2}>0.\]
This is not possible if $u^{*}$ is a smooth solution. $\Box$

For any $m>\frac{4}{3}$ define
\[w_{m}:=1-a_{N,m}r^{\frac{4}{3}}+b_{N,m}r^{m},\]
where
\[a_{N,m}:=\frac{m(N+m-2)}{m(N+m-2)-\frac{4}{3}(N-2/3)}  \ \ \hbox{and} \ \ b_{N,m}:=\frac{\frac{4}{3}(N-2/3)}{m(N+m-2)-\frac{4}{3}(N-2/3)}.\]

Now we are ready to prove the main result of this section.

\begin{theorem}\label{main.sing}  The following upper bounds on $\lambda^*$ hold in large dimensions.
\begin{enumerate}
\item If $N \ge 31$, then Lemma \ref{sing-lem} holds with $u:=w_2$, $\lambda_{N}'=27 \bar \lambda$ and $\beta=\frac{H_N}{2}>27 \bar \lambda$.
\item If $16 \le N \le 30$, then Lemma \ref{sing-lem} holds with $u:=w_3$, $\lambda_{N}'=\frac{H_N}{2}-1, \beta_N=\frac{H_N}{2}$.
\item If $10\le N \le 15$, then Lemma \ref{sing-lem} holds with $u:=w_3$, $\lambda'_N<\beta_N$ given in Table \ref{table:summary}.
\item If $N=9$, then Lemma \ref{sing-lem} holds with $u:=w_{2.8}$, $\lambda'_9:=249<\beta_9:=251$.

\end{enumerate}
The extremal solution is therefore singular for dimensions $N\geq
9$.

\end{theorem}

{\bf Proof.} 1) Assume first that $N\geq 31$, then it is easy to see
that $a_{N,2}<3$ and $a_{N,2}^3\bar{\lambda}\leq27\bar{\lambda}<
\frac{H_{N}}{2}$. We shall show that $w_{2}$ is a singular
$\mathcal{H}-$weak sub-solution of $(P)_{a_{N,2}^3\bar{\lambda}}$
which is stable. Note that $w_{2}\in H^2(B)$, $\frac{1}{1-w_{2}}\in
L^3(B)$, $0\leq w_{2}\leq1$ in $B$, and
\[\Delta^2 w_{2}\leq \frac{a_{N,2}^3\bar{\lambda}}{(1-w_{2})^2} \ \hbox{in} \ \ B \setminus \{0\}.\]
So $w_{2}$ is a $\mathcal{H}-$weak sub-solution of
$(P)_{27\bar{\lambda}}$. Moreover,
\[w_{2}=1-|x|^{\frac{4}{3}}+(a_{N,2}-1)(|x|^{\frac{4}{3}}-|x|^2)\leq 1-|x|^{\frac{4}{3}}.\]
Since $27\bar{\lambda}\leq \frac{H_{N}}{2}$, we get that
\[54\bar{\lambda}\int_{B}\frac{\varphi^2}{(1-w_{2})^3}\leq H_{N}\int_{B}\frac{\varphi^2}{(1-w_{2})^3}\leq H_{N}\int_{B}\frac{\varphi^2}{|x|^4}\leq \int_{B}(\Delta \varphi)^2\]
for all $\varphi \in C^{\infty}_{0}(B)$. Hence, $w_{2}$ is stable.
Thus it follows from Lemma \ref{sing-lem} that $u^{*}$ is singular
and $\lambda^{*}\leq 27\bar{\lambda}$.

2) Assume $16 \leq N\leq 30$ and consider
\[w_{3}:=1-a_{N,3}r^{\frac{4}{3}}+b_{N,3}r^{3}.\]
We show that it is a singular $\mathcal{H}-$weak sub-solution of
$(P_{\frac{H_{N}}{2}-1})$ which is stable. Indeed, we clearly have
$0\leq w_{3}\leq1$ a.e. in $B$, $w_{3} \in H^2(B)$ and
$\frac{1}{1-w_{3}} \in L^3(B)$. Note that
\begin{eqnarray*}
H_{N}\int_{B}\frac{\varphi^2}{(1-w_{3})^3}&=&H_{N}\int_{B}\frac{\varphi^2}{(a_{N,m}r^{\frac{4}{3}}-b_{N,m}r^{m})^3}\\&\leq&
\sup_{0<r<1}
\frac{H_{N}}{(a_{N,m}-b_{N,m}r^{m-\frac{4}{3}})^3}\int_{B}\frac{\varphi^2}{r^4}=H_{N}\int_{B}\frac{\varphi^2}{r^4}\leq
\int_{B}(\Delta \varphi)^2.
\end{eqnarray*}
Using maple one can verify that for $16\leq N\leq31$
\[\Delta^2w_{3}\leq \frac{\frac{H_{N}}{2}-1}{(1-w_{3})^2} \ \ \hbox{on} \ \ (0,1).\]
Hence $w_{3}$ is a sub-solution of $(P_{\frac{H_{N}}{2}-1})$. By
Lemma \ref{sing-lem} $u^{*}$ is singular and $\lambda^{*}\leq
\frac{H_{N}}{2}-1$.

3) Assume $10 \leq N\leq 15$. We shall show that $w_{3}$ satisfies
the assumptions of Lemma \ref{sing-lem} for each dimension $10\leq
N\leq 15$. Using maple, for each dimension $10\leq N\leq 15$, one
can verify that  inequality (\ref{cond1}) holds for $\lambda'_{N}$
given by Table \ref{table:summary}. Then, by using maple again, we
show that there exists $\beta_{N}>\lambda'_{N}$ such that
\begin{eqnarray*}
\frac{(N-2)^2(N-4)^2}{16}\frac{1}{(|x|^2-\frac{N}{2(N-1)}|x|^{\frac{N}{2}+1})(|x|^2-|x|^{\frac{N}{2}})}&+&\frac{(N-1)(N-4)^2}{4}
\frac{1}{|x|^2(|x|^2-|x|^{\frac{N}{2}})}\\ &\geq&
\frac{2\beta_{N}}{(1-w_{3})^3}.
\end{eqnarray*}
The above inequality and improved Hardy-Rellich inequality
(\ref{HR11}) guarantee that the stability condition (\ref{cond2})
holds for $\beta_{N}>\lambda'$. Hence by Lemma \ref{sing-lem} the
extremal solution is singular for $10\leq N\leq 15$. The values of
$\lambda_{N}$ and $\beta_{N}$ are shown in Table
\ref{table:summary}.

\begin{table}[ht]
\caption{Summary} 
\centering 
\begin{tabular}{c c c } 
\hline\hline 
N & $\lambda'_{N}$ & $\beta_{N}$  \\ [0.5ex] 
\hline 
9 & 249 & 251 \\
10 & 320 & 367  \\
11 & 405 & 574  \\
12 & 502 & 851 \\
13 & 610 & 1211 \\
14 & 730 & 1668 \\
15 & 860 & 2235  \\
$16 \leq N \leq 30$ & $\frac{H_{N}}{2}-1$ & $\frac{H_{N}}{2}$\\
$N\geq 31$ & $27 \bar{\lambda}$ & $\frac{H_{N}}{2}$\\ [1ex] 
\hline 
\end{tabular}
\label{table:summary} 
\end{table}

4) Let u:=$w_{2.8}$. Using Maple on can see that
\[\Delta ^2 u \leq \frac{249}{(1-u)^2} \ \ \hbox{in} \ \ B\]
and
\[\frac{502}{(1-u(r))^3}\leq W(r) \ \ \hbox{for all} \ \ r\in(0,1), \]
where $W$ is given by (\ref{W}). Since, $502>2\times 249$, by Lemma
\ref{sing-lem} the extremal solution $u^{*}$ is singular in
dimension $N=9$.  $\Box$
\begin{remark} It follows from the proof of Theorem \ref{main.sing}
that for $N\geq9$ and $\frac{\tau}{\beta}$ sufficiently small, there
exists $u \in H^{2}(B)\cap W^{4,\infty}_{loc}(B\setminus \{0\})$
such that
\begin{equation}\label{cond1}
\Delta^2 u-\frac{\tau}{\beta}\Delta u \leq
\frac{\lambda''_{N}}{(1-u)^2} \ \ \hbox{for} \ \ 0<r<1,
\end{equation}

\begin{equation}
u(1)=0, \ \ \ \ \Delta u|_{r=1}=0,
\end{equation}
\begin{equation}
\hbox{u is singular},
\end{equation}
and
\begin{equation}\label{cond2}
2\beta'_{N} \int_{B}\frac{ \varphi^2}{(1-u)^3}\leq \int_{B}(\Delta
\varphi)^2+\frac{\tau}{\beta}|\nabla \varphi|^2\ \ \hbox{for all} \
\ \varphi \in H^2(B)\cap H^1_0(B),
\end{equation}

where $\beta'_{N}>\lambda''_{N}>0$ are constants. Indeed, for each
dimension $N\geq 9$, it is enough to take $u$ to be the sub-solution
we constructed in the proof of Theorem \ref{main.sing},
$\beta'_{N}:=\beta_{N}$, $\lambda'<\lambda''<\beta$. If
$\frac{\tau}{\beta}$ is sufficiently small so that
$-\frac{\tau}{\beta}\Delta u<\frac{\lambda''-\lambda'}{(1-u)^2}$ on
$(0,1)$, then with an argument similar to that of Lemma
\ref{sing-lem} we deduce that the extremal solution $u^{*}$ of
$(P_{\lambda, \beta, \tau, 0, 0})$ is singular. We believe that the
extremal solution of $(P_{\lambda, \beta, \tau, 0, 0})$ is singular
for all $\beta, \tau>0$ in dimensions $N\geq9$.

\end{remark}

\section{Appendix: Improved Hardy-Rellich Inequalities}

We now prove the improved Hardy-Rellich inequalities used in section
4. They rely on the results of Ghoussoub-Moradifam in \cite{GM}
which provide necessary and sufficient conditions for such
inequalities to hold. At the heart of this characterization is the
following notion of a Bessel pair of functions.
\begin{definition} Assume that $B$ is a ball of radius $R$ in $\R^N$, $V,W \in C^{1}(0,1)$, and
$\int^{R}_{0}\frac{1}{r^{N-1}V(r)}dr=+\infty$. Say that the couple
$(V, W)$  is a {\it Bessel pair on $(0, R)$} if the ordinary
differential equation
\begin{equation*}
\hbox{ $({\rm B}_{V,W})$  \quad \quad \quad \quad \quad \quad \quad
\quad \quad \quad \quad
$y''(r)+(\frac{N-1}{r}+\frac{V_r(r)}{V(r)})y'(r)+\frac{W(r)}{V(r)}y(r)=0$
\quad \quad \quad \quad \quad \quad \quad \quad \quad \quad \quad}
\end{equation*}
has a positive solution on the interval $(0, R)$.
\end{definition}
The needed inequalities will follow from the following two results.
\begin{theorem}\label{GM1} {\bf (Ghoussoub-Moradifam \cite{GM}) } Let $V$ and $W$ be positive radial $C^1$-functions   on $B\backslash \{0\}$, where $B$ is a ball centered at zero with radius $R$ in $\R^N$ ($N \geq 1$) such that  $\int^{R}_{0}\frac{1}{r^{N-1}V(r)}dr=+\infty$ and $\int^{R}_{0}r^{N-1}V(r)dr<+\infty$. The following statements are then equivalent:

\begin{enumerate}

\item $(V, W)$ is a Bessel pair on $(0, R)$.

\item $ \int_{B}V(|x|)|\nabla \phi |^{2}dx \geq \int_{B} W(|x|)\phi^2dx$ for all $\phi \in
C^\infty_{0}(B)$.\\
\end{enumerate}
\end{theorem}

\begin{theorem} \label{main.hr} Let $B$ be the unit ball in $\R^N$ ($N \geq
5$). Then the inequality
\begin{equation}\label{gen-rel}
\hbox{  \quad \quad  $\int_{B}|\Delta u |^{2}dx \geq  \int_{B} \frac{|\nabla u|^{2}}{|x|^2-\frac{N}{2(N-1)}|x|^{\frac{N}{2}+1}} 
dx+(N-1)\int_{B}\frac{|\nabla u|^2}{|x|^2}dx,$  \quad \quad \quad 
\quad \quad \quad \quad  } \\
\end{equation}
holds for all $u \in C^{\infty}_0(\bar{B})$.
\end{theorem}

We shall need the following result to prove (\ref{gen-rel}).\\
\begin{lemma}\label{b.hardy}
For every $u \in C^{1}([0,1])$ the following inequality holds
\begin{equation}\label{bhardy}
\int^{1}_{0}|u'(r)|^{2}r^{N-1} dr \geq \int^{1}_{0}
\frac{u^2}{r^2-\frac{N}{2(N-1)}r^{\frac{N}{2}+1}}
r^{N-1}dr-(N-1)(u(1))^2.
\end{equation}

\end{lemma}
Proof. Let $\varphi:=r^{-\frac{N}{2}+1}-\frac{N}{2(N-1)}$ and
$k(r):=r^{N-1}$. Define $\psi(r)=u(r)/\varphi(r)$, $r \in [0,1]$.
Then
\begin{eqnarray*}
\int^{1}_{0}|u'(r)|^{2}k(r) dr&=&\int^{1}_{0} |\psi (r)|^{2}|\varphi'(r)|^{2}k(r)dr+\int^{1}_{0} 2\varphi(r)\varphi'(r)\psi(r)\psi'(r)k(r)dr +\int^{1}_{0} |\varphi (r)|^2|\psi' (r)|^{2}k(r)dr \\
&=&\int^{1}_{0} |\psi (r)|^{2}(|\varphi'(r)|^{2}k(r)-( k\varphi
\varphi')'(r))dr+\int^{1}_{0} |\varphi (r)|^2|\psi'
(r)|^{2}k(r)dr+\psi^2(1)\varphi'(1)\varphi(1)\\
&\geq&\int^{1}_{0} |\psi (r)|^{2}(|\varphi'(r)|^{2}k(r)-( k\varphi
\varphi')'(r))dr+ \psi^2(1)\varphi'(1)\varphi(1)
\end{eqnarray*}
Note that
$\psi^2(1)\varphi'(1)\varphi(1)=u^2(1)\frac{\varphi'(1)}{\varphi(1)}=-(N-1)u^2(1)$.
Hence, we have
\begin{eqnarray}
\int^{1}_{0}|u'(r)|^{2}k(r) dr&\geq&\int^{1}_{0} -u^2(r)\frac{k'(r)\varphi'(r)+k(r)\varphi''(r)}{\varphi})dr-(N-1)u^2(1)\\
\end{eqnarray}
Simplifying the above inequality we get (\ref{bhardy}). $\square$ \\

The decomposition of a function into its spherical harmonics will be
one of our tools to prove Theorem \ref{main.hr}. Let $u \in
C^{\infty}_0(\bar{B})$. By decomposing $u$ into spherical 
harmonics we get
\[
\hbox{$u=\Sigma^{\infty}_{k=0}u_{k}$ where
$u_{k}=f_{k}(|x|)\varphi_{k}(x)$}
\]
 and $(\varphi_k(x))_k$ are the orthonormal eigenfunctions of the Laplace-Beltrami operator  
with corresponding eigenvalues $c_{k}=k(N+k-2)$, $k\geq 0$. The functions $f_{k}$ belong to 
$u \in C^{\infty}([0,1])$, $f_{k}(1)=0$, and satisfy $f_{k}(r)=O(r^k)$ and $f'(r)=O(r^{k-1})$ as $r \rightarrow 
0$. In particular,
\begin{equation}\label{zero}
\hbox{ $\varphi_{0}=1$ and  $f_{0}=\frac{1}{N
\omega_{N}r^{N-1}}\int_{\partial B_{r}}u ds= \frac{1}{N
\omega_{N}}\int_{|x|=1}u(rx)ds.$}
\end{equation}
We also have  for any $k\geq 0$, and any continuous real valued $W$ on 
$(0,1)$,
\begin{equation}
\int_{B}|\Delta u_{k}|^{2}dx=\int_{B}\big( \Delta 
f_{k}(|x|)-c_{k}\frac{f_{k}(|x|)}{|x|^2}\big)^{2}dx,
\end{equation}
and
\begin{equation}\label{g.formu}
\int_{B}W(|x|)|\nabla u_{k}|^{2}dx=\int_{B}W(|x|)|\nabla 
f_{k}|^{2}dx+c_{k}\int_{B}W(|x|)|x|^{-2}f^{2}_{k}dx.
\end{equation}

Now we are ready to prove Theorem \ref{main.hr}. We shall use the
inequality

\begin{equation}\label{1-dim}
\hbox{$\int^{1}_{0}|x'(r)|^2r^{N-1}dr \geq \frac{(N-2)^2}{4}\int^{1}_{0}\frac{x^{2}(r)}{r^2-\frac{N}{2(N-1)}r^{\frac{N}{2}+1}}r^{N-1}dr$ for all $x\in 
C^1([0,1])$ with $x(1)=0$.}
\end{equation}

{\bf Proof of Theorem \ref{main.hr}:} For all $N\geq 5$ and $k\geq
0$ we have
\begin{eqnarray*}
\frac{1}{Nw_N}\int_{B}|\Delta u_{k}|^{2}dx&=&\frac{1}{Nw_N}\int_{B}\big( \Delta 
f_{k}(|x|)-c_{k}\frac{f_{k}(|x|)}{|x|^2}\big)^{2}dx\\
&=&
\int^{1}_{0}\big(f_{k}''(r)+\frac{N-1}{r}f_{k}'(r)-c_{k}\frac{f_{k}(r)}{r^2}\big)^{2}r^{N-1}dr
\\
&=&\int^{1}_{0}(f_{k}''(r))^{2}r^{N-1}dr+(N-1)^{2} \int^{1}_{0}(f_{k}'(r))^{2}r^{N-3}dr\\
&&+c^{2}_{k} \int^{1}_{0}f_{k}^{2}(r)r^{N-5}
+ 2(N-1) \int^{1}_{0}f_{k}''(r)f_{k}'(r)r^{N-2}\\
&&-2c_{k} \int^{1}_{0}f_{k}''(r)f_{k}(r)r^{N-3}dr - 2c_{k}(N-1)
\int^{1}_{0}f_{k}'(r)f_{k}(r)r^{N-4}dr.
\end{eqnarray*}
Integrate by parts and use (\ref{zero}) for $k=0$ to get
\begin{eqnarray}
\frac{1}{N\omega_{N}}\int_{B}|\Delta u_{k}|^{2}dx&\geq&  
\int^{1}_{0}(f_{k}''(r))^{2}r^{N-1}dr+(N-1+2c_{k}) \int^{1}_{0}(f_{k}'(r))^{2}r^{N-3}dr 
\label{piece.0}\\
&+&
(2c_{k}(n-4)+c^{2}_{k})\int^{1}_{0}r^{n-5}f_{k}^{2}(r)dr+(N-1)(f_{k}^{'}(1))^2\nonumber
\end{eqnarray}
Now define $g_{k}(r)=\frac{f_{k}(r)}{r}$ and note that $g_{k}(r)=O(r^{k-1})$ for all $k\geq 1$. We 
have
\begin{eqnarray*}
\int^{1}_{0}(f_{k}'(r))^{2}r^{N-3}&=&\int^{1}_{0}(g_{k}'(r))^{2}r^{N-1}dr+\int^{1}_{0}2g
_{k}(r)g_{k}'(r)r^{N-2}dr+\int^{1}_{0}g_{k}^{2}(r)r^{N-3}dr\\
&=&\int^{1}_{0}(g_{k}'(r))^{2}r^{N-1}dr-(N-3)\int^{1}_{0}g_{k}^{2}(r)r^{N-3}dr
\end{eqnarray*}
Thus,
\begin{equation}\label{g1}
\int^{1}_{0}(f_{k}'(r))^{2}r^{N-3}\geq 
\frac{(N-2)^2}{4}\int^{1}_{0}\frac{f_{k}^{2}(r)}{r^2-\frac{N}{2(N-1)}r^{\frac{N}{2}+1}}r^{N-3}dr-(N-3)\int^{1}_{0}f_{k}^{2}(r)r^{N-5}dr
\end{equation}
Substituting $2c_k\int^{1}_{0}(f_{k}'(r))^{2}r^{N-3}$ in (\ref{piece.0}) by its lower estimate 
in the last inequality (\ref{g1}), and using Lemma \ref{b.hardy} we
get
\begin{eqnarray*}
\frac{1}{N\omega_{N}}\int_{B}|\Delta u_{k}|^{2}dx&\geq&
\frac{(N-2)^2}{4}\int^{1}_{0}\frac{(f_{k}'(r))^{2}}{r^2-\frac{N}{2(N-1)}r^{\frac{N}{2}+1}}r^{N-1}dr+2c_{k}\frac{(N-2)^2}{4}\int^{1}_{0}\frac{f_{k}^{2}(r)}{r^2-\frac{N}{2(N-1)}r^{\frac{N}{2}+1}}r^{n-3}dr \\  &+&(N-1) 
\int^{1}_{0}(f_{k}'(r))^{2}r^{N-3}dr+c_{k}(N-1) \int^{1}_{0}(f_{k}(r))^{2}r^{N-5}dr\\
&+&c_{k}(c_{k}-(N-1))\int^{1}_{0}r^{N-5}f_{k}^{2}(r)dr+
c_{k}\int^{1}_{0}\frac{(N-2)^2}{4(r^2-\frac{N}{2(N-1)}r^{\frac{N}{2}+1})}-\frac{2}{r^2})dr.\\
&\geq & \frac{(N-2)^2}{4}\int^{1}_{0}\frac{(f_{k}'(r))^{2}}{r^2-\frac{N}{2(N-1)}r^{\frac{N}{2}+1}}r^{N-1}dr+c_{k}\frac{(N-2)^2}{4}\int^{1}_{0}\frac{f_{k}^{2}(r)}{r^2-\frac{N}{2(N-1)}r^{\frac{N}{2}+1}}r^{n-3}dr \\  &+&(N-1) 
\int^{1}_{0}(f_{k}'(r))^{2}r^{N-3}dr+c_{k}(N-1) \int^{1}_{0}(f_{k}(r))^{2}r^{N-5}dr\\
\end{eqnarray*}
The proof is complete in the view of (\ref{g.formu}). \hfill $\square$ \\

We shall now deduce the following corollary.
\begin{corollary}\label{GM2}
Let $N\geq 5$ and $B$ be the unit ball in $\R^N$. Then the following
improved Hardy-Rellich inequality holds for all $\phi \in H^2(B)\cap
H^1_0(B)$:
\begin{eqnarray}\label{HR11}
\int_{B}(\Delta \phi)^2 &\geq&
\frac{(N-2)^2(N-4)^2}{16}\int_{B}\frac{\phi^2}{(|x|^2-\frac{N}{2(N-1)}|x|^{\frac{N}{2}+1})(|x|^2-|x|^{\frac{N}{2}})}\\
&&+\frac{(N-1)(N-4)^2}{4}
\int_{B}\frac{\phi^2}{|x|^2(|x|^2-|x|^{\frac{N}{2}})}\nonumber
\end{eqnarray}
\end{corollary}
{\bf Proof.} Let $\alpha:=\frac{N}{2(N-1)}$ and
$V(r):=\frac{1}{r^2-\alpha r^{\frac{N}{2}+1}}$ and note that
$$\frac{V_{r}}{V}=-\frac{2}{r}+\frac{\alpha(N-2)}{2} \frac{r^{\frac{N}{2}-2}}{1-\alpha r^{\frac{N}{2}-1}}\geq-\frac{2}{r}.$$
The function $y(r)=r^{-\frac{N}{2}+2}-1$ is decreasing and is then a
positive super-solution on $(0,1)$ for the ODE
\begin{equation*}
y''+(\frac{N-1}{r}+\frac{V_{r}}{V})y'(r)+\frac{W_{1}(r)}{V(r)}y=0,
\end{equation*}
where
\[W_{1}(r)=\frac{(N-4)^2}{4(r^2-r^{\frac{N}{2}})(r^2-\alpha r^{\frac{N}{2}+1})}.\]
Hence, by Theorem \ref{GM1} we deduce
$$\int_{B}\frac{|\nabla
\phi|^2}{|x|^2-\alpha|x|^{\frac{N}{2}+1}}\geq(\frac{N-4}{2})^2\int_{B}\frac{\phi^2}{(|x|^2-\alpha
|x|^{\frac{N}{2}+1})(|x|^2-|x|^{\frac{N}{2}})}$$ for all $\phi \in
H^2(B)\cap H^1_0(B)$. Similarly, for $V(r)=\frac{1}{r^2}$ we have
that
$$\int_{B}\frac{|\nabla \phi|^2}{|x|^2}\geq
(\frac{N-4}{2})^2\int_{B}\frac{\phi^2}{|x|^2(|x|^2-|x|^{\frac{N}{2}})}$$
for all $\phi \in H^2(B)\cap H^1_0(B)$. Combining the above two
inequalities with (\ref{bhardy}) we get  (\ref{HR11}). $\Box$

\begin{corollary}\label{GM2} Let $N\geq 7$ and $B$ be the unit ball
in $\R^N$. Then the following improved Hardy-Rellich inequality
holds for all $\phi \in H^2(B)\cap H^1_0(B)$:

\begin{equation}\label{HR3}
\int_{B}|\Delta u|^2\geq \int_{B}W(|x|)u^2.
\end{equation}
where
\begin{equation}\label{W}
W(r)=K(r)(\frac{(N-2)^2}{4(r^2-\frac{N}{2(N-1)}r^{\frac{N}{2}+1})}+
\frac{(N-1)}{r^2}),
\end{equation}
\[K(r)=-\frac{\varphi''(r)+\frac{(n-3)}{r}\varphi'(r)}{\varphi(r)},\]
and
\[\varphi(r)=r^{-\frac{N}{2}+2}+9r^{-2}+10r-20.\]
\end{corollary}
{\bf Proof.} Let $\alpha:=\frac{N}{2(N-1)}$ and
$V(r):=\frac{1}{r^2-\alpha r^{\frac{N}{2}+1}}$. Then $\varphi$ is a
sub-solution for the ODE
\begin{equation*}
y''+(\frac{N-1}{r}+\frac{V_{r}}{V})y'(r)+\frac{W_{2}(r)}{V(r)}y=0,
\end{equation*}
where
\[W_{2}(r)=\frac{K(r)}{r^2-\alpha r^{\frac{N}{2}+1}},\]
Hence by Theorem \ref{GM1} we have
\begin{equation}
\int_{B}\frac{|\nabla u|^2}{|x|^2-\alpha|x|^{\frac{N}{2}+1}}\geq
\int_{B}W_{2}(|x|)u^2.
\end{equation}
Similarly
\begin{equation}
\int_{B}\frac{|\nabla u|^2}{|x|^2}\geq \int_{B}W_{3}(|x|)u^2.
\end{equation}
where
\[W_{3}(r)=\frac{K(r)}{r^2}.\]
Combining the above two inequalities with (\ref{bhardy})  we get
improved Hardy-Rellich inequality (\ref{HR3}). $\Box$ \\

{\bf Acknowledgment:} I would like to thank Professor Nassif
Ghoussoub, my supervisor, for his valuable suggestions, constant
support, and encouragement. I also thank C. Cowan for useful
discussions.

\end{document}